\documentclass[12pt, leqno]{article}
\usepackage{amsmath}
\usepackage{amsthm}
\usepackage{amssymb}
\usepackage{latexsym}
\allowdisplaybreaks[1]
\usepackage{amsfonts}
\usepackage{ascmac}
\usepackage{color}
\usepackage{enumerate}
\usepackage{graphicx}
\usepackage{bm}

\theoremstyle{definition}
\theoremstyle{plain}
\allowdisplaybreaks[1]
\date{}
\newcommand{\p}{\partial}

\newcommand{\dis}{\displaystyle}

\newcommand{\norm}{\parallel}

\newcommand{\N}{{\mathbb N}}
\newcommand{\R}{{\mathbb R}}

\newcommand{\ep}{\varepsilon }

\newcommand{\mmL}{ \mathcal{L}}
\def\text#1{\mbox{#1 }}

\title{\bf A remark on Tonelli's calculus of variations}
\author{Kohei Soga
\footnote{Department of Mathematics, Faculty of Science and Technology, Keio University, 3-14-1 Hiyoshi, Kohoku-ku, Yokohama, 223-8522, Japan. E-mail:  soga@math.keio.ac.jp 
}}
\begin{document}

\maketitle
\begin{abstract} 
\noindent This paper provides a quite simple method of  Tonelli's calculus of variations with positive definite and superlinear Lagrangians. The result complements the classical  literature of calculus of variations before Tonelli's modern approach.  Inspired by Euler's spirit, the proposed method employs finite dimensional approximation of the exact action functional, whose minimizer is easily found as a solution of Euler's discretization of  the exact Euler-Lagrange equation. The Euler-Cauchy polygonal line generated by the approximate minimizer converges to an exact smooth minimizing curve. This framework yields an elementary proof of the existence and regularity of minimizers within the  family of smooth curves and hence, with a minor additional step, within the  family of  Lipschitz curves, without using modern functional analysis on absolutely continuous curves and lower semicontinuity of action functionals.

\medskip\medskip

\medskip


\end{abstract}
%

\setcounter{section}{0}
\setcounter{equation}{0}
\section{Introduction}

Investigation of existence and regularity of minimizing curves of action functionals is one of the most fundamental problems in classical mechanics, which leads to various applications in many fields.  Leonida Tonelli established a general method of the problem, known as the {\it direct method of calculus of variations},  based on the properties of absolutely continuous curves and lower semicontinuity of action functionals. 
We revisit well-established Tonelli's calculus of variations with a Lagrangian $L=L(x,t,\xi)$ satisfying the following conditions:  
\begin{itemize}
\item[(L1)] $L:\R^d\times\R\times \R^d \to \R$, $C^2$;
\item[(L2)] $\frac{\p^2 L}{\p\xi^2}(x,t,\xi)$ is positive definite for each $(x,t,\xi)\in\R^d\times\R\times\R^d$, where  $\frac{\p^2 L}{\p\xi^2}=\Big[\frac{\p^2 L}{\p \xi^i\p\xi^j}\Big]$ stands for the Hessian matrix of $L$ with respect to the $\xi$-variable;
\item[(L3)] $L$ is uniformly superlinear, i.e.,  for each $a\ge0$, there exists $b_a\in\R$ such that $L(x,t,\xi)\ge a|\xi|+b_a$ for all $(x,t,\xi)\in\R^d\times\R\times \R^d$;
\item[(L4)]  The Euler-Lagrange flow $\phi^{t,t^0}_L$ generated by $L$ is complete, i.e,  $\phi^{t,t^0}_L$ is global in time. 
\end{itemize} 
Note that $\phi^{t,t^0}_L$ is well-defined under (L1)-(L3), which is indirectly seen via the corresponding Hamiltonian flow stated in (H1)-(H4) below. 
Let $w:\R^d\to\R$ be a given $C^0$-function such that there exist constants $\alpha\ge0$ and $\beta\in\R$ for which $w(x)\ge -\alpha|x|+\beta$ holds for all $x\in\R^d$.  Let $AC([0,t];\R^d)$ be the family of all absolutely continuous curves $\gamma:[0,t]\to\R^d$. For each $t>0$, we consider the action functional  
\begin{eqnarray*}
\mmL_w:AC([0,t];\R^d)\to \R\cup\{+\infty\},\,\,\,
\mmL_w(\gamma):=\int_0^tL(\gamma(s),s,\gamma'(s))ds+w(\gamma(0)). 
\end{eqnarray*}
We write $\mmL_0$ if $w=0$. Since $L-b_0\ge0$ in $\R^d\times\R\times\R^d$ due to (L3), the integral in $\mmL_w(\gamma)$ always makes sense including its value $+\infty$. 
The typical minimizing problems are 
\begin{eqnarray}\label{AC-minimizing1}
\inf_{\gamma\in AC([0,t];\R^d),\gamma(t)=x}\mmL_w(\gamma),\,\,\,\,\,\,\,\,
\inf_{\gamma\in AC([0,t];\R^d),\gamma(t)=x,\gamma(0)=\tilde{x}}\mmL_0(\gamma), 
\end{eqnarray}
where $x,\tilde{x}\in\R^d$ are arbitrarily given end points. 
The first one in \eqref{AC-minimizing1} is also called the {\it Bolza problem}. 
We refer to \cite{G} for the history and earlier developments of minimization problems for $\mmL_w$. 
Before Tonelli, existence of minimizers was not directly stated, but necessary conditions of minimizers were widely discussed.  
Tonelli's theory states that
\begin{itemize}
\item  {\it Under (L1)--(L4), each of the minimizing problems \eqref{AC-minimizing1} admits at least one minimizing curve in $AC([0,t];\R^d)$;}
\item {\it  Any minimizing curve of \eqref{AC-minimizing1} has $C^2$-regularity and solves on $[0,t]$ the Euler-Lagrange equation generated by $L$, i.e.,
\begin{eqnarray}\label{1EL}
\frac{d}{ds} \frac{\p L}{\p\xi} (\gamma(s),s,\gamma'(s))=\frac{\p L}{\p x}(\gamma(s),s,\gamma'(s)).
\end{eqnarray} 
}  
 \end{itemize}
 We refer to \cite{Clarke} for a good review of Tonelli's theory;  Chapter 3 of \cite{Fathi-book} for full details of the proof of Tonelli's theory;  \cite{Mather},   \cite{Cannarsa} and  \cite{Fathi-book} for important applications of Tonelli's theory to Lagrangian/Hamiltonian dynamics, Hamilton-Jacobi equations and optimal control.    

In actual applications mentioned above, it is often enough to consider  minimization in the smooth class or  in the continuous piecewise smooth class. The Lipschitz class would be also natural, since $\mmL_w$ becomes $\R$-valued. 
Nevertheless, Tonelli's approach requires the absolutely continuous class.  We briefly recall Tonelli's method  together with the reason why the absolutely continuous class is technically necessary, apart from the point that minimization in a wider class of curves would be mathematically more interesting. The first important fact in Tonelli's  approach is that the sequence  $\{ \gamma_j\}_{j\in\N}\subset\{ \gamma\in AC([0,t];\R^d)\,|\, \gamma(t)=x\}$ such that 
\begin{eqnarray*}
\mmL_w(\gamma_j)\to \inf_{\gamma\in AC([0,t];\R^d),\gamma(t)=x}\mmL_w(\gamma)\mbox{\quad as $j\to\infty$}
\end{eqnarray*}
is uniformly bounded and equicontinuous, providing a uniformly convergent subsequence and its limit $\gamma^\ast\in C^0([0,t];\R^d)$. The next important fact is that  the sequence  $\{ \gamma'_j\}_{j\in\N}$ is uniformly integrable, which implies that $\gamma^\ast$ belongs to $AC([0,t];\R^d)$ and $\gamma'_j$ converges to $\gamma^\ast{}'$ in the weak-$L^1$ topology.  Then, the lower semicontinuity of $\mmL_w$ with respect to the $C^0$-topology  resulting from (L2) shows 
$$\mmL_w(\gamma^\ast)\le\liminf_{j\to\infty}\mmL_w(\gamma_j),$$
 to conclude that $\gamma^\ast$ is a minimizing curve of  \eqref{AC-minimizing1}.  
In this argument, it is not a priori clear if $\gamma^\ast$ is Lipschitz or not, even if we start with the Lipschitz class instead of $AC([0,t];\R^d)$. In a subset of the Lipschitz class with an artificial bound of Lipschitz constants, Tonelli's approach may find a Lipschitz minimizing curve with the same artificial bound. However, variations around the minimizer cannot necessarily be taken in all directions and the Euler-Lagrange equation fails to be derived (one does not a priori know if a Lipschitz constant of the minimizer is strictly smaller than the artificial bound).  Tonelli's approach based on compactness of a minimizing sequence and the lower semicontinuity is powerful enough to extend the existence result to the case with Lagrangians out of (L1)--(L4): see, e.g., \cite{BGH} for more details. 

 Higher regularity of $\gamma^\ast$ is easily seen, if there are additional growth conditions of the derivatives $\p_x L$, $\p_\xi L$ that guarantee integrability of $\p_x L(\gamma(s),s,\gamma'(s))$ and $\p_\xi L(\gamma(s),s,\gamma'(s))$ over $[0,t]$ for $\gamma\in AC([0,t];\R^d)$. In this case, we easily obtain the vanishing  G\^ateaux derivative of $\mmL_w$ at $\gamma^\ast$, which gives the Euler-Lagrange equation in the integral form. This concludes that $\gamma^\ast$ has $C^2$-regularity and satisfies the Euler-Lagrange equation. 
 Only with the assumptions (L1)--(L4), one can still obtain the same regularity result through a classical result  by Weierstrass: if $t>0$ and $|x-\tilde{x}|$ are sufficiently small,  the second problem of \eqref{AC-minimizing1} admits the unique $C^2$-minimizer satisfying the Euler-Lagrange equation. See Appendix of \cite{Mather} and Chapter 3 of \cite{Fathi-book} for more details.      

In this paper, we propose a completely different approach to find minimizing curves of $\mmL_w$ directly within the smooth class, which immediately leads to results on the minimization within the Lipschitz class. This is an attempt to complement the missing part of the classical calculus of variations (i.e., existence of minimizers) only by means of classical techniques. {\it We emphasize that  the proposed method does not reduce the significance of Tonelli's approach, also because it works only with Lagrangians satisfying (L1)--(L4). However, the proof of existence and regularity of minimizers for \eqref{AC-minimizing1}  becomes much simpler. Furthermore, it would be interesting to observe that Euler's classical idea to find minimizers (see Chapter 2 of \cite{G}) seems almost successful and general enough. }          

The key idea of the proposed method is to introduce a finite dimensional problem corresponding to  \eqref{AC-minimizing1}, i.e., we discretize $[0,t]$ and deal with the Riemann sum of the integral in $\mmL_w$. The finite dimensional problem is merely minimization of a smooth function defined in $\R^{d'}$ with $d'\in\N$ becoming arbitrarily large for the limit of approximation but finite in each step. Then, we easily obtain a minimizing point with  the vanishing  derivative.
 Furthermore, the vanishing  derivative yields a finite difference equation that is the Euler's discretization of the exact Euler-Lagrange equation \eqref{1EL}; the minimizing point solves the discrete Euler-Lagrange equation with a priori boundedness. 
 Due to (L1)--(L3), we have the Legendre transform $H$ of $L$. The discrete Euler-Lagrange equation is transformed to the  Euler's discretization of exact Hamilton's equations generated by $H$, i.e., 
\begin{eqnarray}\label{1HS} 
 \gamma'(s)=\frac{\p H}{\p p}(\gamma(s),s,p(s)),\quad 
  p'(s)=-\frac{\p H}{\p x}(\gamma(s),s,p(s)).
 \end{eqnarray}
It is straightforward to see that the Euler-Cauchy polygonal line generated by the solution of discrete Hamilton's equations converges to a $C^2$-curve $(\gamma^\ast(s),p(s))$ solving \eqref{1HS}; $\gamma^\ast(s)$ is a $C^2$-minimizer of $\mmL_w$ within the smooth class solving \eqref{1EL}. The essential point of  finite dimensional approximation  is that one can immediately obtain the approximate Euler-Lagrange equation with nice estimates on the discrete level without being bothered by regularity issues. 
 
In order to extend the minimization to the Lipschitz class, we use the standard fact that each Lipschitz curve is approximated by a series of smooth curves with bounded derivatives. Then, Lebesgue's dominated convergence theorem concludes that a minimizing curve within the smooth class is also minimizing among all Lipschitz curves. It is not difficult to see that any other minimizing curve in the Lipschitz class (if it exists) necessarily has $C^2$-regularity and solves the Euler-Lagrange equation.   
 Note that in the case of  absolutely continuous curves,  derivatives are unbounded in general and Lebesgue's dominated convergence theorem would not work; hence our approach without any additional growth condition of $L$ would fail in $AC([0,t];\R^d)$.  


\setcounter{section}{1}
\setcounter{equation}{0}
\section{Minimization in the smooth class}

We discuss minimization of $\mmL_w$ in the smooth class as simply as possible. We proceed in a slightly wider class than $C^1([0,t];\R^d)$.  Let $C^1_{\rm pw}([0,t];\R^d)$ be  the family of all continuous piecewise smooth curves; namely, $C^1_{\rm pw}([0,t];\R^d)$ is the family of all continuous curves $\gamma:[0,t]\to\R^d$ such that each $\gamma$ has a finite division $0=\tau_0<\tau_1<\cdots<\tau_m=t$ with $m=m(\gamma)\ge 1$ for which $\gamma$ is $C^1$ within each interval $[\tau_{l-1},\tau_l]$. 
Note that each $\gamma\in C^1_{\rm pw}([0,t];\R^d)$ is Lipschitz. 
Calculus of variations for $\mmL_w$ in $C^1_{\rm pw}([0,t];\R^d)$ is the minimum requirement in analysis of viscosity solutions (or the value functions) to the Hamilton-Jacobi equation generated by $H$ and weak KAM theory based on the Lax-Oleinik type operator generated by $L$. Reasoning in $C^1_{\rm pw}([0,t];\R^d)$ is essentially the same as that in $C^1([0,t];\R^d)$ (in both cases, one does not need the technique of mollification and Lebesgue's dominated convergence theorem). Note that our reasoning works also in $C^2([0,t];\R^d)$.   
In this section, we discuss how to find a $C^2$-minimizing curve directly within $C^1_{\rm pw}([0,t];\R^d)$. The result of Section 3 shows that any minimizer of $C^1_{\rm pw}([0,t];\R^d)$ necessarily has $C^2$-regularity and solves the Euler-Lagrange equation.   

We briefly refer to the Legendre transform $H$ of $L$ with respect to the $\xi$-variable:  
\begin{eqnarray}\label{Legandre}
H(x,t,p)=\sup_{\xi\in\R^d}\{p\cdot\xi-L(x,t,\xi)\}.
\end{eqnarray}
Here, $x\cdot y=\sum_{i=1}^dx^iy^i$ for $x,y,\in\R^d$. It follows from (L1)--(L4) that   
\begin{enumerate}
\item[(H1)] $H:\R^d\times\R\times\R^d\to\R$, $C^2$;
\item[(H2)] $\frac{\p^2H}{\p p^2}(x,t,p)$ is positive definite for each $(x,t,p)\in\R^d\times\R\times\R^d$;
\item[(H3)] $H$ is locally uniformly superlinear with respect to $p$, i.e., for each $a\ge0$ and $K\subset\R^d\times\R$ compact, there exists $b'_{a,K}\in\R$ such that $H(x,t,p)\ge a|p|+b'_{a,K}$  for all $(x,t,p)\in K\times \R^d$;
\item[(H4)] The Hamiltonian flow $\phi^{t,t_0}_H$ generated by $H$ is complete.
\end{enumerate} 
Furthermore, since $\p_\xi L(x,t,\xi)=p$ is invertible with respect to the $\xi$-variable for any $p\in\R^d$ due to (L1)--(L3), we have the $C^1$-map $\xi=\xi(x,t,p)$ such that $\p_\xi L(x,t,\xi(x,t,p))\equiv p$, where the supremum of \eqref{Legandre} is uniquely attained by $\xi=\xi(x,t,p)$.  This implies that  
\begin{eqnarray*}
&&\frac{\p H}{\p p}(x,t,\cdot):\R^d\to\R^d\,\,\,\mbox{is bijection},\,\,\,\,\,\frac{\p H}{\p p}(x,t,p)=\xi(x,t,p),\\
&&\frac{\p H}{\p p}(x,t,p)=\xi\Leftrightarrow p=\frac{\p L}{\p \xi}(x,t,\xi),\,\,\,\,\,
\frac{\p H}{\p x}(x,t,p)=-\frac{\p L}{\p x}(x,t,\xi(x,t,p)),\\
&&\mbox{$\gamma(s)$ solves \eqref{1EL}$\Leftrightarrow$ the pair $\gamma(s)$, $\dis p(s)=\frac{\p L}{\p\xi}(\gamma(s),s,\gamma'(s))$ solves \eqref{1HS}}.
\end{eqnarray*}
Note that the uniqueness property of the Euler-Lagrange equation \eqref{1EL} is derived from that of Hamilton's equations, while the completeness of the Hamiltonian flow follows from (L4); if $\sup\{L(x,t,\xi)\,|\,x\in\R^d,t\in\R, |\xi|\le r \}<+\infty$ for each $r\ge0$, then $H$ is uniformly superlinear and $\sup\{H(x,t,\xi)\,|\,x\in\R^d,t\in\R, |\xi|\le r \}<+\infty$ for each $r\ge0$.   
We refer to \cite{Cannarsa} and \cite{Fathi-book} for detailed description of the Legendre transform and to Lemma 3.1 of \cite{Soga5} for more details of the above statements.    
\subsection{Problem with one fixed end point and one free end point}
For each $t>0$ and $x\in\R^d$, consider the minimizing problem with one fixed end point and one free end point:
\begin{eqnarray}\label{C^1-minimizing}
\inf_{\gamma\in C_{\rm pw}^1([0,t];\R^d),\gamma(t)=x}\mmL_w(\gamma).
\end{eqnarray}

\medskip
\noindent{\bf Claim.} {\it One can find a $C^2$-minimizing curve that attains \eqref{C^1-minimizing} and satisfies  the Euler-Lagrange equation generated by $L$, without going through the argument in $AC([0,t];\R^d)$ and Weierstrass's theorem. }
\medskip
\medskip

\noindent Our claim is justified by reasoning with  finite dimensional approximation of $\mmL_w$:  
\begin{eqnarray*}
&&h:=\frac{t}{K}, \,\,\, t_k:=hk,\,\,\, K,k\in\N\,\,\,\mbox{ \quad (we will send $K\to\infty$)},\\
&&\bm{y}=(y_0,y_1,\ldots,y_{K-1})\in \R^{dK},\,\,\,\,\,\,y_k\in\R^d\mbox{ for $k=0,1,\ldots,K-1$},\\
&& \bm{y}'=(y'_0,y'_1,\ldots,y'_{K-1}),\,\,\,y'_k:=\frac{y_{k+1}-y_{k}}{h}\mbox{ for $k=0,1,\ldots,K-1$ with } \,\,\,y_K=x,\\
&&\mmL_w^K:\R^{dK}\to\R,\,\,\,\mmL_w^K(\bm{y}):=\sum_{k=0}^{K-1}L(y_k,t_k,y'_k)h+w(y_0).
\end{eqnarray*}
\indent {\bf Step 1.} We minimize $\mmL_w^K$ in $\R^{dK}$. 
It follows from (L3) and the condition $w(x)\ge -\alpha|x|+\beta$ that, for $a=\alpha$ and any $\bm{y}\in\R^{dK}$, we have the lower bound of $\mmL_w^K$ as 
\begin{eqnarray*}
\mmL_w^K(\bm{y})&\ge& \sum_{k=0}^{K-1}(\alpha |y'_k|+b_\alpha) h-\alpha|y_0|+\beta
\ge \alpha\Big|\sum_{k=0}^{K-1} y'_k\Big|h +b_\alpha t-\alpha|y_0|+\beta\\
&=&\alpha|x-y_0|+b_\alpha t-\alpha|y_0|+\beta
\ge -\alpha|x|+b_\alpha t+\beta>-\infty.
\end{eqnarray*}
For $\bm{y}_x:=(x,\ldots,x)\in\R^{dK}$, we have  
\begin{eqnarray}\label{CCCC}
C_x:=\sup_{K\in\N}\mmL_w^K(\bm{y}_x)=\sup_{K\in\N}\sum_{k=0}^{K-1}L(x,t_k,0)h+w(x)<+\infty.
\end{eqnarray}
It is enough to minimize $\mmL_w^K$ within the set 
$Y_x^K:=\{ \bm{y}\in\R^{dK}\,|\,\mmL_w^k(\bm{y})\le C_x \}.$ 
It follows from (L3) that, for $a=1+\alpha$ and any $\bm{y} \in Y_x^K$, we have 
\begin{eqnarray}\label{bound1}
\,\,\,C_x&\ge&\mmL_w^K(\bm{y})\ge \sum_{k=0}^{K-1}|y'_k|h-\alpha|x|+b_{1+\alpha}t+\beta\ge  \sum_{k=l\ge0}^{K-1}|y'_k|h-\alpha|x|+b_{1+\alpha}t+\beta\\\nonumber
&\ge& |y_l|-|x|-\alpha|x|+b_{1+\alpha}t+\beta,\\
\label{bound2}
\,\,\,C_x &\ge&\mmL_w^K(\bm{y})\ge \min_{0\le k\le K-1}|y'_k|t-\alpha|x|+b_{1+\alpha}t+\beta. 
\end{eqnarray}
Hence, we obtain by \eqref{bound1},
\begin{eqnarray}\label{bound3}
|y_l| \le R_1:=C_x+(1+\alpha)|x|-b_{1+\alpha}t-\beta,\,\,\,\forall\,l=0,\ldots,K-1,  \,\,\,\forall\,\bm{y}\in Y_x^K. 
\end{eqnarray}
Hereafter, $R_1,R_2,\ldots$ stand for some constants  independent from $K$.  
 Therefore, $Y_x^K$ is a bounded subset of $\R^{dK}$ and there exists $\bar{\bm{y}}=\bar{\bm{y}}(K)\in Y_x^K$ such that 
$$\mmL_w^K(\bar{\bm{y}})=\inf_{\bm{y}\in\R^{dK}}\mmL_w^K(\bm{y}).$$
It follows from \eqref{bound2} that we have $k^\ast=k^\ast(K)$ such that $|\bar{y}'_{k^\ast(K)}|= \min_{0\le k\le K-1}|\bar{y}'_k|$ and 
\begin{eqnarray}\label{bound4}
|\bar{y}{}_{k^\ast(K)}'|\le R_2:=(C_x+\alpha|x|-b_{1+\alpha}t-\beta)t^{-1}.
\end{eqnarray}

{\bf Step 2.}  We take variations around $\bar{\bm{y}}$ and send $K\to\infty$. Since $\mmL_w^K(y_0,y_1,\ldots,y_{K-1})$ is a $C^1$-function except for the  $y_0$ variable ($w$ is only continuous), its minimum point $\bar{\bm{y}}$ implies 
\begin{eqnarray}\label{ELK}
\frac{\p}{\p y_k} \mmL_w^K(\bar{\bm{y}})&=&\frac{\p L}{\p x}(\bar{y}_k,t_k,\bar{y}{}_k')h-\Big(\frac{\p L}{\p \xi}(\bar{y}_k,t_k,\bar{y}{}_k')-\frac{\p L}{\p \xi}(\bar{y}_{k-1},t_{k-1},\bar{y}{}'_{k-1})\Big)\\\nonumber
&=&0,\,\,\,\forall\,k=1,\ldots,K-1.
\end{eqnarray}
The Legendre transform equivalently transform \eqref{ELK} into \eqref{HSK}:
\begin{eqnarray}\nonumber
&&\bar{z}_k:=\frac{\p L}{\p \xi}(\bar{y}_k,t_k,\bar{y}{}_k')\,\,\,\mbox{ for }k=0,\ldots,K-1,\\\nonumber
&&\bar{z}{}'_k:=\frac{\bar{z}_k-\bar{z}_{k-1}}{h}\,\,\,\,\mbox{ for } k=1,\ldots,K-1,
\\\label{HSK}
&&\left\{
\begin{array}{l}\dis 
\bar{y}{}_k'=\frac{\p H}{\p p}(\bar{y}_k,t_k,\bar{z}_k),\,\,\,\,\,\,\forall\,k=0,\ldots,K-1,\medskip\\\dis
\bar{z}{}'_k=-\frac{\p H}{\p x}(\bar{y}_k,t_k,\bar{z}_k),\,\,\,\,\,\,\forall\,k=1,\ldots,K-1.
\end{array}
\right. 
\end{eqnarray}
Note that \eqref{HSK} is discretization of Hamilton's equations in terms of the Euler method.   Due to \eqref{bound3} and \eqref{bound4}, we have for any $K\in\N$, 
$$|\bar{y}_{k^\ast(K)}|\le R_1,\,\,\,t_{k^\ast(K)}\in[0,t],\,\,\, |\bar{z}_{k^\ast(K)}|\le R_3:=\sup_{|x|\le R_1, 0\le s\le t,|\xi|\le R_2} \Big|\frac{\p L}{\p\xi}(x,s,\xi)\Big |.$$
Note that, if $w\in C^1$, we can take variations also with respect to $y_0$ and obtain an instant a priori bound of $\bar{z}_0$.  
Since $\{(\bar{y}_{k^\ast(K)},t_{k^\ast(K)},\bar{z}_{k^\ast(K)})\}_{K\in\N}$ is a bounded sequence of $\R^d\times\R\times\R^d$, we find a convergent subsequence (denoted with $\{K_j\}_{j\in\N}$) with its limit $(x^\ast, t^\ast,p^\ast)$ satisfying $|x^\ast|\le R_1,\,\,\,t^\ast\in[0,t],\,\,\, |p^\ast|\le R_3$.   
We assert that the Euler-Cauchy polygon method and (H4) yield the following fact: 
 {\it the polygonal line generated by the solution $\{(\bar{y}_{jk},\bar{z}_{jk})\}_{k=0,\ldots,K_j-1}$ of  \eqref{HSK} with 
$K=K_j$ converges uniformly to 
$$(\gamma^\ast(s),p^\ast(s)):=\phi_H^{s,t^\ast}(x^\ast,p^\ast)$$
 on $[0,t]$ as $j\to\infty$, where $\gamma^\ast$ satisfies the Euler-Lagrange equation with $\gamma^\ast(t)=x$.} 

We check our assertion. It follows from (H4) that  $(\gamma^\ast(s),p^\ast(s))$ exists on the interval $[0,t]$. 
Due to the continuous dependency with respect to initial data, the curve 
$$(\gamma_j^\ast(s),p_j^\ast(s)):=\phi_H^{s,t_{k^\ast(K_j)}}(\bar{y}_{k^\ast(K_j)},\bar{z}_{k^\ast(K_j)})$$
 uniformly converges to $(\gamma^\ast(s),p^\ast(s))$ on $[0,t]$ as $j\to\infty$. 
Hence, $\{(\gamma_j^\ast(s),s,p_j^\ast(s))\}_{s\in[0,t]}$ is contained in a bounded neighborhood $B\subset\R^d\times\R\times\R^d$ of $\{(\gamma^\ast(s),s,p^\ast(s))\}_{s\in [0,t]} $ for all large $j$. 
The Taylor approximation yields with  $h_j:=t/K_j$, 
\begin{eqnarray*}
&&\gamma_j^\ast(t_{k+1})-\gamma_j^\ast(t_{k})=\frac{\p H}{\p p}(\gamma_j^\ast(t_{k}),t_k,p_j^\ast(t_k))h_j +O(h_j^2),\,\,\,k=0,\ldots,K_j-1,\\
&&p^\ast_j(t_{k+1})-p_j^\ast(t_{k})=-\frac{\p H}{\p x}(\gamma_j^\ast(t_{k+1}),t_{k+1},p_j^\ast(t_{k+1}))h_j +O(h_j^2),\,\,\,k=0,\ldots,K_j-1,
\end{eqnarray*}
where $|O(h_j^2)|$ is bounded by $R_4h_j^2$ for any $k$ and $j\to\infty$ with some constant $R_4\ge0$ determined by the supremums of the first and second derivatives of $H$ in $B$.  
Let $\bar{y}_j{}_k,\bar{z}_j{}_k$ be the solution of \eqref{HSK} with $K=K_j$ and set 
$$U_j{}_k=(u_j{}_k,\tilde{u}_j{}_k):=(\bar{y}_j{}_k-\gamma_j^\ast(t_{k}),\bar{z}_j{}_k-p_j^\ast(t_{k})),$$ where $U_{jk^\ast(K_j)}=0$. 
Let $\theta>0$ be a common  Lipschitz constant $\p_p H$ and $\p_x H$ in $B$. As long as $(\bar{y}_{jk},t_k,\bar{z}_{jk})$ belongs to $B$ (this is a priori the case at least for $k$ close to $k^\ast(K_j)$ for all large $j$), we have 
\begin{eqnarray*}
|u_j{}_{k+1}-u_j{}_k|\le \theta|U_j{}_k|h_j+R_4h_j^2,\,\,\,\,\,\,|\tilde{u}_j{}_{k+1}-\tilde{u}_j{}_k|\le \theta|U_j{}_{k+1}|h_j+R_4h_j^2,
\end{eqnarray*}
which lead to 
\begin{eqnarray*}
&&|U_j{}_{k+1}|\le \frac{1+\theta h_j}{1-\theta h_j}|U_{j}{}_k|+R_5h_j^2\le (1+3\theta h_j)|U_{j}{}_k|+R_5h_j^2,\,\,\,\,k\ge k^\ast(K_j),\\
&&|U_j{}_{k-1}|\le \frac{1+\theta h_j}{1-\theta h_j}|U_{j}{}_{k}|+R_5h_j^2\le (1+3\theta h_j)|U_{j}{}_k|+R_5h_j^2,\,\,\,\,k\le k^\ast(K_j),\\
&&\Big(|U_{jk+1}|+\frac{R_5h_j}{3\theta}\Big)\le (1+3\theta h_j)\Big(|U_{jk}|+\frac{R_5h_j}{3\theta}\Big),\,\,\,k\ge k^\ast(K_j),\\
&&\Big(|U_{jk-1}|+\frac{R_5h_j}{3\theta}\Big)\le (1+3\theta h_j)\Big(|U_{jk}|+\frac{R_5h_j}{3\theta}\Big),\,\,\,k\le k^\ast(K_j),\\
&&\Big(|U_{jk}|+\frac{R_5h_j}{3\theta}\Big)\le (1+3\theta h_j)^{k-k^\ast(K^j)}\Big(|U_{jk^\ast(K_j)}|+\frac{R_5h_j}{3\theta}\Big),\,\,\,k\ge k^\ast(K_j),\\
&&\Big(|U_{jk}|+\frac{R_5h_j}{3\theta}\Big)\le (1+3\theta h_j)^{k^\ast(K_j)-k}\Big(|U_{jk^\ast(K_j)}|+\frac{R_5h_j}{3\theta}\Big),\,\,\,k\le k^\ast(K_j).
\end{eqnarray*}
It follows from $(1+\nu)^{\frac{1}{\nu}}\to e$ as $\nu\to0+$ (the convergence is from below) that
\begin{eqnarray*}
(1+3\theta h_j)^{|k-k^\ast(K^j)|}\le (1+3\theta h_j)^{K_j}= (1+3\theta h_j)^{\frac{t}{h_j}}=\Big\{(1+3\theta h_j)^{\frac{1}{3\theta h_j}}\Big\}^{3\theta t}\le e^{3\theta t}.
\end{eqnarray*}
Therefore,  we find  a constant $R_6>0$ independent of $j$ for which  
\begin{eqnarray*}
|U_{j}{}_k|\le R_6 h_j,\,\,\,\mbox{ $\forall\,k=k_1,k_1+1,\ldots, k_2$},  
\end{eqnarray*}
where $k_1,k_2$ are such that $0\le k_1\le k^\ast(K_j)\le k_2\le K_j-1$ and $\{(\bar{y}_{jk},t_k,\bar{z}_{jk})\}_{k_1\le k\le k_2}\subset B$. 
Now it is clear that $\{(\bar{y}_{jk},t_k,\bar{z}_{jk})\}_{0\le k\le K_j-1}\subset B$ for all large $j$ and the Euler-Cauchy polygonal line generated by $\{(\bar{y}_{jk},\bar{z}_{jk})\}_{0\le k\le K_j-1}$  uniformly converges  to $(\gamma^\ast(\cdot),p^\ast(\cdot))$ as $j\to\infty$. 
In particular, we have 
\begin{eqnarray*}
\max_{0\le k\le K_j-1}|\bar{y}_j{}'_{k}- \gamma^\ast{}'(t_k)  |&=&\max_{0\le k\le K_j-1}\Big|\frac{\p H}{\p p}(\bar{y}_{jk},t_k,\bar{z}_{jk})- \frac{\p H}{\p p}(\gamma^\ast(t_k),t_k,p^\ast(t_k)) \Big|\\
&\to& 0\quad \mbox{as $j\to\infty$}. 
\end{eqnarray*}
Thus, denoting $\bar{\bm{y}}_j:=(\bar{y}_j{}_{0},\bar{y}_j{}_{1},\ldots,\bar{y}_j{}_{K_j-1})$, we obtain 
\begin{eqnarray*}
\mmL_w^{K_j}(\bar{\bm{y}}_j)\to \mmL_w(\gamma^\ast)\quad \mbox{as $j\to\infty$}. 
\end{eqnarray*} 

{\bf Step 3.} We check that $\gamma^\ast$ attains \eqref{C^1-minimizing}.  For each $\gamma\in C^1_{\rm pw}([0,t];\R^d)$ with $\gamma(t)=x$, there exists a division $0=\tau_0<\tau_1<\cdots<\tau_{m-1}<\tau_m=t$ with  $m=m(\gamma)\ge 1$ for which $\gamma$ is $C^1$ on each interval $[\tau_{l-1},\tau_l]$. Let $K_j$ and $h_j$ be the one given in Step 2. For each $K_j$, define $\bm{y}_j:=(\gamma(t_0),\gamma(t_1),\ldots, \gamma(t_{K_j-1}))$. Note that each component  $y_j'{}_k=(\gamma(t_{k+1})-\gamma(t_k))/h_j$ of $\bm{y}_j'$ satisfies for any $k$ with  a quantity $\eta(r)\ge0$ tending to $0$ as $r\to0+$, 
\begin{eqnarray*}
|y'_j{}_k|\le R_7:=\sup_{s\in[0,t]}|\gamma'(s)|,\,\,\,\mbox{ if $\exists\,\tau_l\in(t_k,t_{k+1})$;\,\,\, otherwise \,\,\,\,}|y'_j{}_k-\gamma'(t_k)|\le \eta(h_j). 
\end{eqnarray*}
 Let $\Lambda_j$ be the set of $k\in\{0,1,\ldots, K_j-1\}$ such that the interval $(t_k,t_{k+1})$ does not contain any $\tau_l$. 
We have
\begin{eqnarray*}
|\mmL_w(\gamma)-\mmL_w^{K_j}(\bm{y}_j)|&\le& 
\sum_{k\in\Lambda_j}\int_{t_k}^{t_{k+1}}|L(\gamma(s),s,\gamma'(s))-L(y_{jk},t_k,y'_{jk})|ds\\
&&+R_8(m-1)h_j\to0\mbox{\quad as $j\to\infty$}, \\
\mmL_w(\gamma)-\mmL_w(\gamma^\ast)&=&
\mmL_w(\gamma)- \mmL_w^{K_j}(\bm{y}_j)+\mmL_w^{K_j}(\bm{y}_j)-\mmL_w^{K_j}(\bar{\bm{y}}_j)+\mmL_w^{K_j}(\bar{\bm{y}}_j)-  \mmL_w(\gamma^\ast)\\
&\ge&
\mmL_w(\gamma)- \mmL_w^{K_j}(\bm{y}_j)+\mmL_w^{K_j}(\bar{\bm{y}}_j)-  \mmL_w(\gamma^\ast)\to 0\quad \mbox{ as $j\to\infty$}, 
\end{eqnarray*}
which implies $\mmL_w(\gamma^\ast)\le \mmL_w(\gamma)$. 
\subsection{Problem with two fixed end points}
For each $t>0$ and $x,\tilde{x}\in\R^d$, consider the minimizing problem with two fixed end points:
\begin{eqnarray}\label{2-C^1-minimizing}
\inf_{\gamma\in C_{\rm pw}^1([0,t];\R^d),\gamma(t)=x, \gamma(0)=\tilde{x}}\mmL_0(\gamma).
\end{eqnarray}
 We may obtain a minimizing curve of \eqref{2-C^1-minimizing} in the same way as the previous subsection, where we fix $y_0=\tilde{x}$ and take $\bm{y}=(y_1,\ldots,y_{K-1})\in\R^{d(K-1)}$ as free variables and use 
 $$\bm{y}_{x\tilde{x}}:=(\gamma(t_1),\ldots,\gamma(t_{K-1}))\mbox{ with }\gamma(s):=\tilde{x}+\frac{x-\tilde{x}}{t}s$$
in \eqref{CCCC}, instead of $\bm{y}_x=(x,\ldots,x)\in\R^{dK}$.  
 
\setcounter{section}{2}
\setcounter{equation}{0}
\section{Minimization in the Lipschitz class}
Let ${\rm Lip}([0,t];\R^d)$ be the family of all Lipschitz continuous curves. 
Let  $\gamma^\ast$ be a minimizing curve of \eqref{C^1-minimizing}. 
We show that any $\gamma\in {\rm Lip}([0,t];\R^d)$ with $\gamma(t)=x$ satisfies 
$\mmL_w(\gamma^\ast)\le \mmL_w(\gamma)$. 
Let  $M>0$ be a constant such that $|\gamma'(s)|\le M$ a.e. $s\in[0,t]$.  We extend $\gamma$ to be $x$ for $s>t$ and to be $\gamma(0)$ for $s<0$, and   
mollify it with the standard mollifier with a small parameter $\ep>0$, obtaining a smooth curve $\tilde{\gamma}^\ep:[0,t]\to\R^d$ with $\tilde{\gamma}^\ep{}':[0,t]\to [-M,M]^d$ for each $\ep>0$. We set 
$$\gamma^\ep(s):= \tilde{\gamma}^\ep(s)+ (x-\tilde{\gamma}^\ep(t))\frac{s}{t}+(\gamma(0)-\tilde{\gamma}^\ep(0))\frac{t-s}{t}:[0,t]\to\R^d.$$
Note  that $\gamma^\ep(t)=x$,  $\gamma^\ep(0)=\gamma(0)$ for each $\ep>0$ and 
$\norm \gamma^\ep-\gamma\norm_{C^0([0,t];\R^d)}\to0$,  $\norm \gamma^\ep{}'-\gamma'\norm_{L^1([0,t];\R^d)}\to0$ as $\ep\to0+$. 
Since $\sup_{s\in[0,t]}|L(\gamma^\ep(s),s,\gamma^\ep{}'(s))|$ is bounded independently from $\ep\to0+$,   Lebesgue's dominated convergence theorem yields    
$\lim_{\ep\to0+}\mmL_w(\gamma^\ep)=\mmL_w(\gamma)
$, where we take a subsequence of $\ep\to0+$ so that  $\gamma^\ep{}'(s)\to\gamma'(s)$ a.e. $s\in[0,t]$ (this holds even when $w$ is lower semicontinuous). 
Since $\gamma^\ep$ belongs to $C^1_{\rm pw}([0,t];\R^d)$ with $\gamma^\ep(t)=x$, we have
$\mmL_w(\gamma^\ast)\le\mmL_w(\gamma^\ep),
$
which yields our assertion. 

\indent Finally, we check that any other possible minimizer of $\mmL_w$ in $ {\rm Lip}([0,t];\R^d)|_{\gamma(t)=x}$ has $C^2$-regularity and satisfies  the Euler-Lagrange equation generated by $L$.  
Let $\gamma\in {\rm Lip}([0,t];\R^d)$ be  such that  $\gamma(t)=x$ and $\mmL_w(\gamma)=\mmL_w(\gamma^\ast)$. Then, for any smooth curve $v:[0,t]\to\R^d$ with $v(t)=v(0)=0$, we have 
\begin{eqnarray*}
0&=&\frac{d}{d\ep}\mmL_w(\gamma+\ep v)|_{\ep=0}\\
&=& \int_0^t \Big\{\frac{\p L}{\p\xi}(\gamma(s),s,\gamma'(s)) -\int_0^s\frac{\p L}{\p x}(\gamma(\tau),\tau,\gamma'(\tau))d\tau\Big\}\cdot v'(s)\,ds,
\end{eqnarray*} 
which implies 
$$\frac{\p L}{\p\xi}(\gamma(s),s,\gamma'(s)) -\int_0^s\frac{\p L}{\p x}(\gamma(\tau),\tau,\gamma'(\tau))d\tau=\mbox{constant,\quad  a.e. $s\in[0,t]$}.$$
Since $\p_\xi L(x,s,\cdot)$ is invertible and $\int_0^s\p_x L(\gamma(\tau),\tau,\gamma'(\tau))d\tau$ is continuous as a function of $s\in[0,t]$, we see that $\gamma'(s)$ is continuous. Then, $\int_0^s\p_x L(\gamma(\tau),\tau,\gamma'(\tau))d\tau$ turns into  $C^1$ to conclude that $\gamma$ has $C^2$-regularity and satisfies the Euler-Lagrange equation.

We may deal with the problem with two fixed end points in the same way, where we extend $\gamma\in {\rm Lip}([0,t];\R^d)$ to be $x$ for $s>t$ and to be $\tilde{x}$ for $s<0$, and  take 
$$\gamma^\ep(s):= \tilde{\gamma}^\ep(s)+ (x-\tilde{\gamma}^\ep(t))\frac{s}{t}+(\tilde{x}-\tilde{\gamma}^\ep(0))\frac{t-s}{t}:[0,t]\to\R^d.$$

\medskip\medskip

\noindent{\bf Acknowledgement.} The author  is supported by JSPS Grant-in-aid for Young Scientists \#18K13443.

\end{document}